\RequirePackage{fix-cm}
	\documentclass[twocolumn]{svjour3}          
\smartqed  
\usepackage{graphicx}
\usepackage{amssymb}
\usepackage{amsmath,color,latexsym} 

\usepackage{extarrows,verbatim}
\usepackage[normalem]{ulem}

 \journalname{}

\begin{document}

\title{Structurally stable families of periodic solutions in sweeping processes of networks of elastoplastic springs 
}


\author{Ivan Gudoshnikov, Oleg Makarenkov 
}


\institute{           O. Makarenkov \at
              Department of Mathematical Sciences, University of Texas at Dallas, 75080 Richardson, USA\\
              \email{makarenkov@utdallas.edu}
}

\date{Received: date / Accepted: date}

\maketitle

\begin{abstract} Networks of elastoplastic springs (elastoplastic systems) have been linked to differential equations with polyhedral constraints in the pioneering paper by Moreau (1974). Periodic loading of an elastoplastic system, therefore, corresponds to a periodic motion of the polyhedral constraint. According to Krejci (1996), every solution of a sweeping process with a periodically moving constraint asymptotically converges to a periodic orbit.  Understanding whether such an asymptotic periodic orbit is unique or there can be an entire family of asymptotic periodic orbits (that form a periodic attractor) has been an open problem since then. Since suitable small perturbation of a polyhedral constraint seems to be always capable to destroy a potential family of periodic orbits, it is expected that none of potential periodic attractor is  structurally stable. In the present paper we give a simple example to prove that even though the periodic attractor (of non-stationary periodic solutions) can be destroyed by little perturbation of the moving constraint, the periodic attractor resists perturbations of the physical parameters of the mechanical model (i.e. the parameters of the network of elastoplastic springs).  


\end{abstract}

\section{Introduction} Networks of elastoplastic springs are increasingly used in the modeling of the distribution of stresses in elastopastic media \cite{buxton,chen}, swarming of mobile router networks \cite{robot2,robot1}, and  other physical phenomena. According to Moreau \cite{moreau}, the stresses of springs of such a network can be described by a differential inclusion
({\it Moreau sweeping process})  
\begin{equation}
-y'(t)\in N_{C(t)}(y(t)),\quad y(t)\in\mathbb{R}^m,
\label{eq:MSP0}
\end{equation}
where $C(t)\subset\mathbb{R}^m$ is a closed polyhedron that plays the role of a constraint,
\begin{equation}\nonumber
  N_{{C}}(x)=\left\{\begin{array}{ll}\left\{\zeta \in \mathbb{R}^n:\langle\zeta,c-x\rangle\leqslant 0,\  c\in {{C}}\right\},& {\rm if}\ x\in {{C}},\\
   \emptyset,& {\rm if}\ x\not\in {{C}},
\end{array}\right.
\end{equation}
and the dimension $m$ equals or smaller than the number of springs in the network.

\vskip0.2cm

\begin{figure}[h]
\begin{center}  
\includegraphics[scale=0.78]{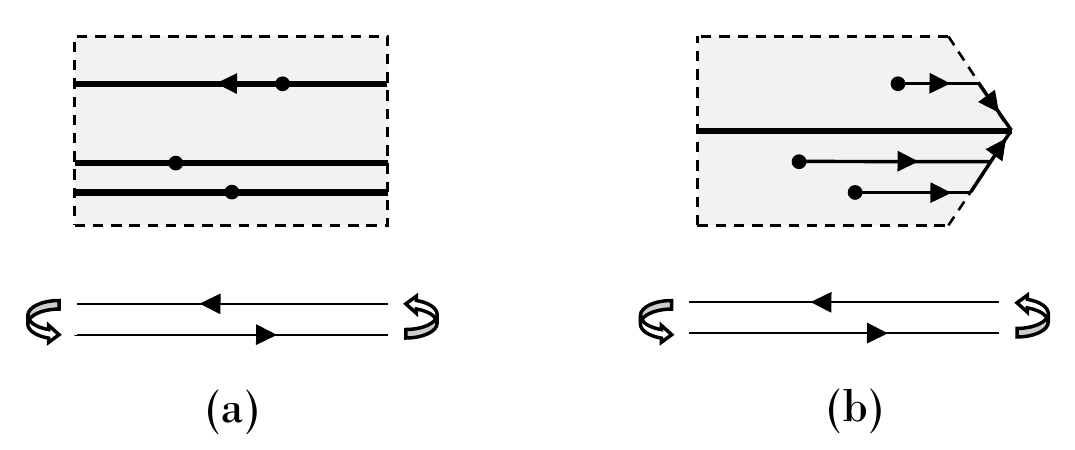}  
\end{center}\vskip0.0cm 
\caption{Sample trajectories (solid curves) of Moreau sweeping process with a moving constraint (dashed rectangle) that moves back and forth. The bold points are the initial conditions of  sample trajectories. The figure illustrates the type of attractor (solid black curves) when (a) the moving constraint is just a rectangle, (b) the moving constraint is a pentagon with a corner that accumulates all the trajectories.}
\label{fig1}
\end{figure}

\noindent Periodicity of the constraint $C(t)$ corresponds to periodicity of the external loading applied to the given network of springs. The fundamental result by Krejci \cite[Theorem 3.14]{Krejci1996} says that for $C(t)$ of the form $C(t)=C+c(t)$, where $C$ is a convex closed bounded set and $t\mapsto c(t)$ is a $T$-periodic vector-function, any solution of sweeping process (\ref{eq:MSP0}) converges to some $T$-periodic regime. For a class of continuum elastoplastic media with $T$-periodic loading the uniqueness of $T$-periodic response is established in Frederick-Armstrong \cite[p.~159]{frederick}. Sufficient conditions for the uniqueness of the response in sweeping processes can be drawn based on Adly et al \cite{adly}.
The non-uniqueness of the response for sweeping processes can of course be easily designed, see Fig.~\ref{fig1}a, where one gets a family of periodic solutions by moving a rectangle normal to its sides back and worth. However, as shown at Fig.~\ref{fig1}b, small perturbation of such a rectangle destroys the attracting family of 
orbits of Fig.~\ref{fig1}a leaving only a single attracting solution. 


\noindent That is why a natural question arises: 
\begin{itemize}
\item[] whether or not any network of elastoplastic springs can always be slightly perturbed in way that destroys any potential family of periodic orbits in the respective sweeping process (\ref{eq:MSP0})?
\end{itemize}
\begin{figure}[h]\center
\vskip-0.7cm
\includegraphics[scale=0.78]{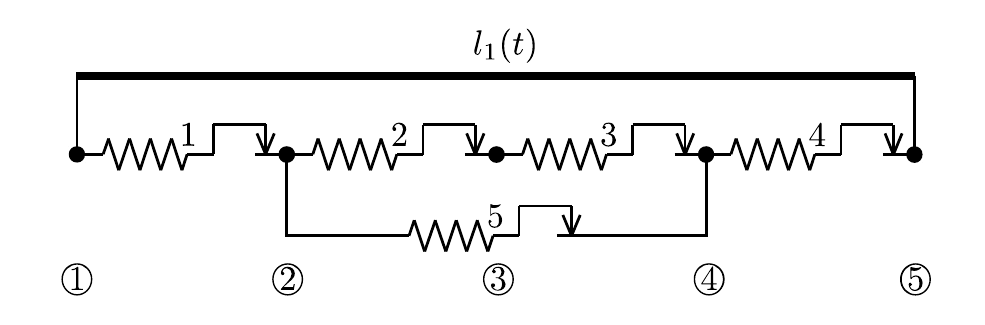}\vskip-0.0cm
\caption{\footnotesize A one-dimensional network of 5 springs on 5 nodes with one displacement-controlled loading. The circled digits stand for numbers of nodes. The regular digits are the numbers of springs. The thick bar is the displacement-controlled loading $l_1(t)$. The stress-controlled loadings $f_1(t),...,f_5(t)$ are applied at nodes.
} \label{ex2newfig}
\end{figure}

\noindent As uniqueness of the response lies in the core of reliability of modeling prediction (see e.g. \cite{bouby,brazil}), the above-stated question is not of merely academic value. We introduce a simple example that answers this question negatively. Specifically, we show that the  cyclically loaded network of elastoplastic springs of Fig.~\ref{ex2newfig}  leads to a sweeping process with a family of attracting periodic orbits. 

\vskip0.2cm

\noindent The paper is organized as follows. In the next section we define a network of elastoplastic springs formally. In section~\ref{eq:sweepProcess} we derive a sweeping process \eqref{eq:MSP0} that governs the quasi-static evolution of such a network. 
Section~\ref{step-by-step} is based on Moreau \cite{moreau} and Gudoshnikov-Makarenkov \cite{G-M}. It compiles a guide for closed-form computation of the quantities required for construction of a sweeping process of a given network of elastoplastic springs. This guide is then used in Section~\ref{sect:structStab} to construct the sweeping process of the network of elastoplastic springs of Fig.~\ref{ex2newfig}. We rigorously proof (Proposition~\ref{prop1} and Corollary~\ref{corr1}) that such a sweeping process admits a family of periodic orbits that persists under perturbations of the mechanical parameters of the network.

\section{A concise definition of a general network of elastoplastic springs}\label{defined}

\begin{figure}[h]\center
\vskip-0.3cm
\includegraphics{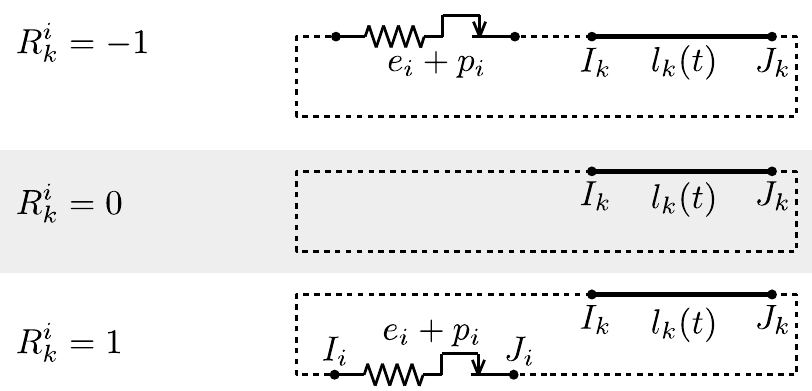}\vskip-0.0cm
\caption{\footnotesize Illustration of the signs of the components of the incidence vector $R^k\in\mathbb{R}^m.$ The dotted contour stays for the chain of the springs associated with the vector $R^k$. 
} \label{figR}
\end{figure}

\noindent We consider a network of $m$ elastoplastic springs on $n$ nodes that are connected according to a directed graph given by the $n\times m$ incidence matrix $D^\top$. The Hooke's coefficients  $a_1, ..., a_m$ of the springs are arranged into an $m\times m$-matrix $A={\rm diag}\left\{a_1,...,a_m\right\}.$ The elastic limits $[c_i^-,c_i^+]$ of springs are used to introduce a parallelepiped $C\subset\mathbb{R}^m$ as $C=[c_1^-,c_1^+]\times...\times[c_m^-,c_m^+].$ In addition the network comes with a collection of stress-controlled and displacement-controlled loadings $\{f_i(t)\}_{i=1}^n$ and $\{l_i(t)\}_{i=1}^q$ respectively. The stress-controlled loadings are simply applied at the $n$ nodes of the network and are supposed to satisfies the equation of static balance
\begin{equation}\label{balance!}
   f_1(t)+...+f_n(t)=0.
\end{equation}
As for the displacement-controlled loading $l_k(t),$ $k\in\overline{1,q}$, we consider a chain of springs which connects the left node $I_k$ of the constraint $k$ with its right node $J_k.$
To each displacement-controlled loading $l_k(t)$ we, therefore, associate a so-called {\it incidence vector}  $R^k\in\mathbb{R}^m$ whose $i$-th component $R^k_i$ is $-1,$ $0,$ or $1$ according to whether the spring $i$ increases, not influences, or decreases the displacement when moving from node $I_k$ to $J_k$ along the chain selected, see Fig.~\ref{figR}. 
We assume that the displacement-controlled loadings $\{l_i(t)\}_{i=1}^m$ are independent in the sense that  
\begin{equation}\label{A2}
{\rm rank}\left(D^\top R\right)=q.
\end{equation}
Mechanically, condition (\ref{A2}) ensures that the displacement-controlled loadings don't contradict one another. For example, (\ref{A2}) rules out the situation where two different displacement-controlled loadings connect same pair of nodes.

\section{A concise formulation of the sweeping process of a general network of elastoplastic springs}
\label{eq:sweepProcess}

\noindent In this section we follow Moreau \cite{moreau} (see also  Gudoshnikov-Makarenkov \cite{G-M}). If condition (\ref{balance!}) holds, then there exists a function 
 $\bar h:\mathbb{R}\to\mathbb{R}^m$, such that
\begin{equation}\label{fh}
   f(t)=-D^\top \bar h(t).
\end{equation}
Then, under condition (\ref{A2}), there exists an $n\times q-$matrix $L$, such that 
\begin{equation}
   R^\top D L=I_{q\times q}.
   \label{barxi}
\end{equation}

\sloppy

\noindent Introducing 
\begin{equation}\label{U}   U=\left\{x\in D\mathbb{R}^n:R^\top x=0\right\},\qquad  V=A^{-1}U^\perp,
\end{equation}
where $U^\perp=\left\{y\in\mathbb{R}^m:\left<x,y\right>=0,\ x\in U\right\},$
the space $V$
becomes an orthogonal complement of the space $U$ in the sense of the scalar product
\begin{equation}\label{product}
   (u,v)_A=\left<u,Av\right>.
\end{equation}
Therefore, any element $x\in\mathbb{R}^m$ can be uniquely decomposed 
 as
$$
   x=P_U x+P_V x,
$$
where $P_U$ and $P_V$ are linear (orthogonal in sense of \eqref{product}) projection maps on $U$ and $V$ respectively.
Define
\begin{align}
  g(t)&=P_V DL l(t),\label{g(t)}\\
 h(t)&=P_U A^{-1}\bar h(t),\label{h}\\ 
  N_C^A(x)&=\\
  &\hskip-1cm =\left\{\begin{array}{ll}\left\{\xi\in\mathbb{R}^m:\left<\xi,A(c-x)\right>\le 0,\  c\in C\right\},&\quad {\rm if}\ x\in C,\\
   \emptyset,& \quad {\rm if}\ x\not\in C,
\end{array}\right.\nonumber\\
  \Pi(t)&=A^{-1}C+h(t)-g(t),\label{Pi(t)}
\end{align}
Assuming that both $f:\mathbb{R}\to\mathbb{R}^n$ and $l:\mathbb{R}\to\mathbb{R}^q$ are Lipschitz continuous, we get that $h(t)$ and $g(t)$ are Lipschitz continuous as well, so that the function 
$$
  y(t)=A^{-1}s(t)+h(t)-g(t)
$$ 
is absolutely continuous for any absolutely continuous $t\mapsto s(t).$

\begin{theorem}\label{thm1} {\rm \cite{moreau} (see also  \cite{G-M})} Assume that the network of elastoplastic springs  $(D,A,C,R,f(t),l(t))$ of section~\ref{defined} satisfies  the conditions (\ref{balance!}) and (\ref{A2}). Assume  that 
$h:\mathbb{R}\to\mathbb{R}^m$ and $g:\mathbb{R}\to\mathbb{R}^m$ given by 
(\ref{g(t)})-(\ref{h}) are Lipschitz continuous. Assume that safe load condition
\begin{equation}\label{A4}
\left(C+Ah(t)\right)\cap U^\perp\not=\emptyset
\end{equation}
holds on some time interval $[0,T].$
Then, 
the function $s(t)=(s_1(t),...,s_m(t))$ defines the evolution of stresses of 
the network  $(D,A,C,R,f(t),l(t))$ for $t\in[0,T]$ if and only if the function 
$$
  y(t)=A^{-1}s(t)+h(t)-g(t)
$$ 
satisfies the differential inclusion (called {\it sweeping process})
\begin{eqnarray}
   -\dot y&\in& N^A_{\Pi(t)\cap V}(y),\quad {\rm for\ a.a.\ }t\in[0,T], \label{sw1}\\
y(0)&\in & \Pi(0)\cap V.\label{eq:initialElastic}
\end{eqnarray}
\end{theorem}

\noindent It remains to note that,  for Lipschitz continuous
$h:\mathbb{R}\to\mathbb{R}^m$ and $g:\mathbb{R}\to\mathbb{R}^m$, 
sweeping process (\ref{balance!})-(\ref{A2}) has a unique Lipschitz-continuous solution for any initial condition (see e.g. Kunze and Monteiro Marques \cite[sect. 3]{kunze}).


\section{The shakedown condition}
\label{sect:conditions}


\noindent {The following} conditions will rule out the existence of constant solutions.

\begin{proposition}\label{prop3} {\rm \cite[Proposition~3]{G-M}} Assume that conditions of Theorem~\ref{thm1} hold. If 
\begin{equation}\label{formulaprop3}
\|A^{-1}c^--A^{-1}c^+\|_A<\left\|g(t_1)-g(t_2)\right\|_A,
\end{equation}
for some $0\le t_1< t_2$, 
where 
\begin{eqnarray*}
   \|x\|_A&=&\sqrt{\left<x,Ax\right>},\\ 
   c^-&=&(c_1^-,...,c_m^-)^\top,\\ 
   c^+&=&(c_1^+,...,c_m^+)^\top.
\end{eqnarray*}
then sweeping process (\ref{sw1}) doesn't have any solutions that are constant on $[t_1,t_2].$
\end{proposition}

\begin{remark}\label{remarkprop3} Note, the left-hand-side in the squared inequality \eqref{formulaprop3} from the statement of Proposition~\ref{prop3} can be computed as \newline
$\|A^{-1}c^--A^{-1}c^+\|_A^2=
 \left<c^--c^+,A^{-1}(c^--c^+)\right>.$
\end{remark}


\section{A step-by-step guide to compute the quantities of the sweeping process from a network of elastoplastic springs}\label{step-by-step}

\noindent In this section we again follow Moreau \cite{moreau}, but use the notations and additional properties established in Gudoshnikov-Makarenkov \cite{G-M}. In particular, \cite[Lemma~1]{G-M} and \cite[formula~(49)]{G-M} say that  
\begin{eqnarray}
\dim U&=&n-q-1.\label{dimU}\\
\dim V&=&m-n+q+1,\label{dimV}
\end{eqnarray}
provided that (\ref{A2}) is satisfied.

\vskip0.2cm

\noindent {\boldmath{\bf Step 1. The matrix $M.$}} According to (\ref{dimU}), there should exist an $n\times (n-q-1)-$matrix $M$ such that 
\begin{equation}\label{RTDM}
  R^\top DM=0 \quad {\rm and}\quad {\rm rank}(DM)=n-q-1
\end{equation}
which allows to introduce $U_{basis}$ as
\begin{equation}\label{Ubasis}
  U_{basis}=DM.
\end{equation} 

\noindent {\boldmath{\bf Step 2. The matrix $V_{basis}.$}}
According to (\ref{U}), $V_{basis}$ is an arbitrary matrix of $m-n+q+1=\dim V$ linearly independent columns that solves
\begin{equation}\label{Vbasis}
  (U_{basis})^\top AV_{basis}=0.
\end{equation}

\noindent {\boldmath{\bf Step 3. The matrix $D^\perp.$}} Define $D^\perp$  to be an $m\times(m-n+1)-$matrix of full rank that solves the equation
\begin{equation}\label{Dperp}
  (D^\perp)^\top D={\color{black}0_{(m-n+1)\times(m-n+1)}}.
\end{equation}

\noindent {\bf Step 4. Other quantities.} Using Steps 2 and 3, we can compute an $(m-n+q+1)\times q$-matrix $\bar L$ as
\begin{equation}\label{barLnew}
   \bar L=\left(\left(\begin{array}{c}
      R^\top \\
      (D^\perp)^\top\end{array}\right)V_{basis}\right)^{-1}\left(\begin{array}{c}
      I_{q\times q}\\ 0_{(m-n+1)\times q}\end{array}\right).
\end{equation}
It turns out that formula (\ref{g(t)}) can now be rewritten in closed-form as
\begin{equation}\label{g(t)ex}
   g(t)=V_{basis}\bar L l(t).
\end{equation}
To account for all possible functions $h(t)$ from (\ref{h}) we will simply take $h(t)$ as
\begin{equation}\label{h(t)}
  h(t)=U_{basis}H(t),
\end{equation}
where $H(t)$ is an arbitrary Lipschitz continuous control input. It is possible to compute $H(t)$ in terms of $f(t)$, but it is not of added value here. 

\vskip0.2cm

\noindent Finally, for $\Pi(t)\cap V$ we have
\begin{equation}\label{222}
   \Pi(t)\cap V=\bigcap_{i=1}^{m}V_i(t),
\end{equation}

\noindent where 
\begin{equation}\label{normals}
\begin{array}{rcl}
V_i(t)&=&\left\{x\in V:c_i^-\hskip-0.05cm+\hskip-0.05cma_ih_i(t)\hskip-0.05cm\le\right.\\
&&\hskip1cm \left.\le\hskip-0.05cm\left<n_i,Ax\hskip-0.05cm+\hskip-0.05cmAg(t)\right>\hskip-0.05cm\le c_i^+\hskip-0.05cm+\hskip-0.05cma_ih_i(t)\right\},\\
 n_i&=&V_{basis}\bar n_i,\\
  n_i&=&\left(\left(\begin{array}{c} R^\top\\ (D^\perp)^\top\end{array}\right)V_{basis}\right)^{-1}\left(\begin{array}{c} R^\top\\ (D^\perp)^\top\end{array}\right)e_i,
\end{array}
\end{equation}
and $e_i\in\mathbb{R}^m$ is the vector with 1 in the $i$-th component and zeros elsewhere.

\section{The sweeping process of the network of elastoplastic springs of Figure~\ref{ex2newfig}}
\label{sect:structStab}

The network of elastoplastic springs of Fig.~\ref{ex2newfig} is given by  
\begin{equation}\label{Dex2new}
D\xi=
\left(\begin{array}{ccccc}
-1 &  1 &  0 &  0 & 0\\
 0 & -1 &  1 &  0 & 0\\
 0 &  0 & -1 &  1 & 0\\
 0 &  0 &  0 & -1 & 1\\
 0 & -1 &  0 &  1 & 0
\end{array}\right)
\left(\begin{array}{c}\xi_1\\
\xi_2\\
\xi_3\\
\xi_4\\
\xi_5
\end{array}\right), \quad R=\left(\begin{array}{l} 1 \\ 1 \\ 1 \\ 1 \\ 0\end{array}\right),
\end{equation}

\noindent some $5\times 5$ diagonal matrix $A$ of Hooke's coefficients and some intervals $[c_i^-,c_i^+]$, $i\in\overline{1,5}$, of elasticity bounds.

\vskip0.2cm

\noindent Formula (\ref{dimU}) leads to
$$
  \dim U=5-1-1=3.
$$
The $5\times 3-$matrix $M$ that solves (\ref{RTDM}) and the respective $5\times 3-$matrix (\ref{Ubasis}) are found as
\begin{equation}\label{UbasisEx2new}
  M=\left(\begin{array}{ccc}
   0 & 0 & 0\\ 
   1 & 0 & 0\\
   0 & 1 & 0\\
   0 & 0 & 1\\
   0 & 0 & 0
   \end{array}\right),
   \quad U_{basis}=\left(\begin{array}{ccc} 
    1 &  0 &  0\\
   -1 &  1 &  0\\
    0 & -1 &  1\\
    0 &  0 & -1\\
   -1 &  0 &  1
   \end{array}\right),
\end{equation}
and 
 $H(t)$ in (\ref{h(t)}) is an arbitrary Lipschitz continuous function from $[0,T]$ to $\mathbb{R}^3.$ 

\vskip0.2cm

\noindent 
According to (\ref{dimV}) and (\ref{Vbasis}), one gets
\begin{equation}\label{VbasisEx2}
\begin{array}{rcl}
   \dim V&=&5-5+1+1=2,\\ 
   V_{basis}&=&\left(\begin{array}{ccc}
1/a_1 & \ &0\\
1/a_2 &\ &-1/a_2\\
1/a_3 &\ &-1/a_3\\
1/a_4 &\ &     0\\ 
    0 &\ & 1/a_5
\end{array}\right)
\end{array}
\end{equation}

\noindent Following Step~3 of section~\ref{step-by-step}, we compute $\dim D^\perp=5-5+1=1$ and the $5\times 1$-dimensional solution of (\ref{Dperp}) is
{\color{black}\begin{equation}\label{0110-1}
   D^\perp=\left(\begin{array}{c}
   0\\ 1\\ 1\\ 0\\ -1
   \end{array}\right).
\end{equation}}
Therefore, according to formula (\ref{barLnew}),  the $2\times1-$matrix $\bar L$ computes as 
\begin{equation}
  \bar L=\left(\left(\begin{array}{ccccc}
1 & 1 & 1 & 1 &  0\\
0 & 1 & 1 & 0 & -1\\
\end{array}\right)V_{basis}\right)^{-1}\left(\begin{array}{c}
1\\0\end{array}\right)
\label{eq:LBarEx2new}
\end{equation}
and by (\ref{g(t)ex}) we get
\begin{equation}\label{relation1}
   g(t)=V_{basis}\left(\left(\begin{array}{ccccc}
1 & 1 & 1 & 1 &  0\\
0 & 1 & 1 & 0 & -1\\
\end{array}\right)V_{basis}\right)^{-1}\left(\begin{array}{c}
1\\0\end{array}\right)l(t).
\end{equation}
On the other hand, formula (\ref{normals}) says that for each $i\in\overline{1,5}$, the normal vector $n_i$ is given by $n_i=$
\begin{equation}\label{relation2}
   V_{basis}\left(\left(\begin{array}{ccccc}
1 & 1 & 1 & 1 &  0\\
0 & 1 & 1 & 0 & -1\\
\end{array}\right)V_{basis}\right)^{-1}\hskip-0.2cm\left(\begin{array}{ccccc}
1 & 1 & 1 & 1 &  0\\
0 & 1 & 1 & 0 & -1\\
\end{array}\right)e_i.
\end{equation}
Note, formulas (\ref{relation1}) and (\ref{relation2}) hold for any $a_i,$ $i\in\overline{1,5}$ and any $c_i^-,c_i^+$, $i\in\overline{1,5}.$ Therefore, we see from formulas (\ref{relation1}) and (\ref{relation2}) that $n_1\parallel g(t)$ and $n_4\parallel g(t)$ for any values of the physical parameters of the network of Fig.~\ref{ex2newfig}. However, at this point we don't know whether or not the normals $n_1$ and $n_4$ have anything to do with the sides of the shape $\Pi(t)\cap V$ given by (\ref{222}), as it may happen that the constraints of \eqref{222} provided by $n_1$ and $n_4$ become redundant for a particular $h(t)$.
\begin{proposition}\label{prop1}  There is an open set of the parameters $a_i,$  $c_i^-,c_i^+$, $i\in\overline{1,5}$, and an open set of Lipschitz-continuous functions $H:[0,T]\mapsto\mathbb{R}^3,$ for which the vectors $n_1$ and $n_4$ are the normal vectors of the two opposite sides of the shape $\Pi(t)\cap V$. In particular, this open set of the parameters contains the point
\begin{equation}\label{para}
\begin{array}{l}
   c_i^-=-1,\ c_i^+=1,\ a_i=1,\\
    H(t)\equiv(-0.5,-0.8,-1)^\top.
\end{array}
\end{equation}
Here $[0,T]$ is an arbitrary chosen domain of the functions $t\mapsto H(t).$
\end{proposition}

\noindent {\bf Proof.} Without loss of generality we can consider $g(t)\equiv 0$. Indeed, since $g(t)$ acts along $V$, $g(t)$ simply translates $\Pi(t)\cap V$  within $V$, so that $g(t)$ doesn't change the shape of $\Pi(t)\cap V$.

\vskip0.2cm

\noindent Plugging (\ref{para}) into  (\ref{UbasisEx2new}) and using (\ref{h(t)}) we get
$$
  h(t)\equiv (-0.5,-0.3,-0.2,1,-0.5)^\top.
$$
Therefore, for the parameters (\ref{para}), formula (\ref{222}) says that $x\in\Pi(t)\cap V$ if and only if
\begin{equation}\label{system}
  \left\{\begin{array}{lllll}
     -1-0.5&\le&\left<n_1,x\right>&\le& 1-0.5,\\
    -1-0.3&\le&\left<n_2,x\right>&\le& 1-0.3,\\
    -1-0.2&\le&\left<n_3,x\right>&\le& 1-0.2,\\
    -1+1&\le&\left<n_4,x\right>&\le& 1+1,\\
    -1-0.5&\le&\left<n_5,x\right>&\le& 1-0.5,
  \end{array}\right.
\end{equation}
where
\begin{equation}\label{n1}  \renewcommand\arraystretch{2} 
\begin{array}{l}n_1=n_4=\dfrac{1}{d}\left(-3,-1,-1,-3,-2\right)^\top,\\ 
n_2=n_3=\dfrac{1}{d}\left(-1,-3,-3,-1,2\right)^\top,\\ 
n_5=\dfrac{1}{d}\left(-2,2,2,-2,-4\right)^\top,\\ 
d=-8.
\end{array}
\end{equation}
Based on (\ref{relation2}), $n_1=n_4$ and $n_2=n_3$. Therefore, 1st and 4th lines of system (\ref{system}) as well as 2nd and 3rd lines combine, that reduces the number of double-sided inequalities to 3. Substituting the expressions (\ref{relation2}) with parameters (\ref{para}) into (\ref{system}) and plugging $x=V_{basis}v,$ where $v\in\mathbb{R}^2,$ system (\ref{system}) reduces to the following system
\begin{equation}\label{system2}
  \left.\begin{array}{lrlcll}
{\rm normals}\  n_1\ {\rm and}\  n_4: \quad&    0&\le&v_1&\le& 0.5,\\
{\rm normals}\  n_2\ {\rm and}\  n_3: \quad&     -1.2&\le&v_1-v_2&\le& 0.7,\\
{\rm normal}\  n_5: &     -1.5&\le&v_2&\le& 0.5.    
  \end{array}\right.
\end{equation}
\begin{figure}[h]\center
\vskip-0.5cm
\includegraphics[scale=0.6]{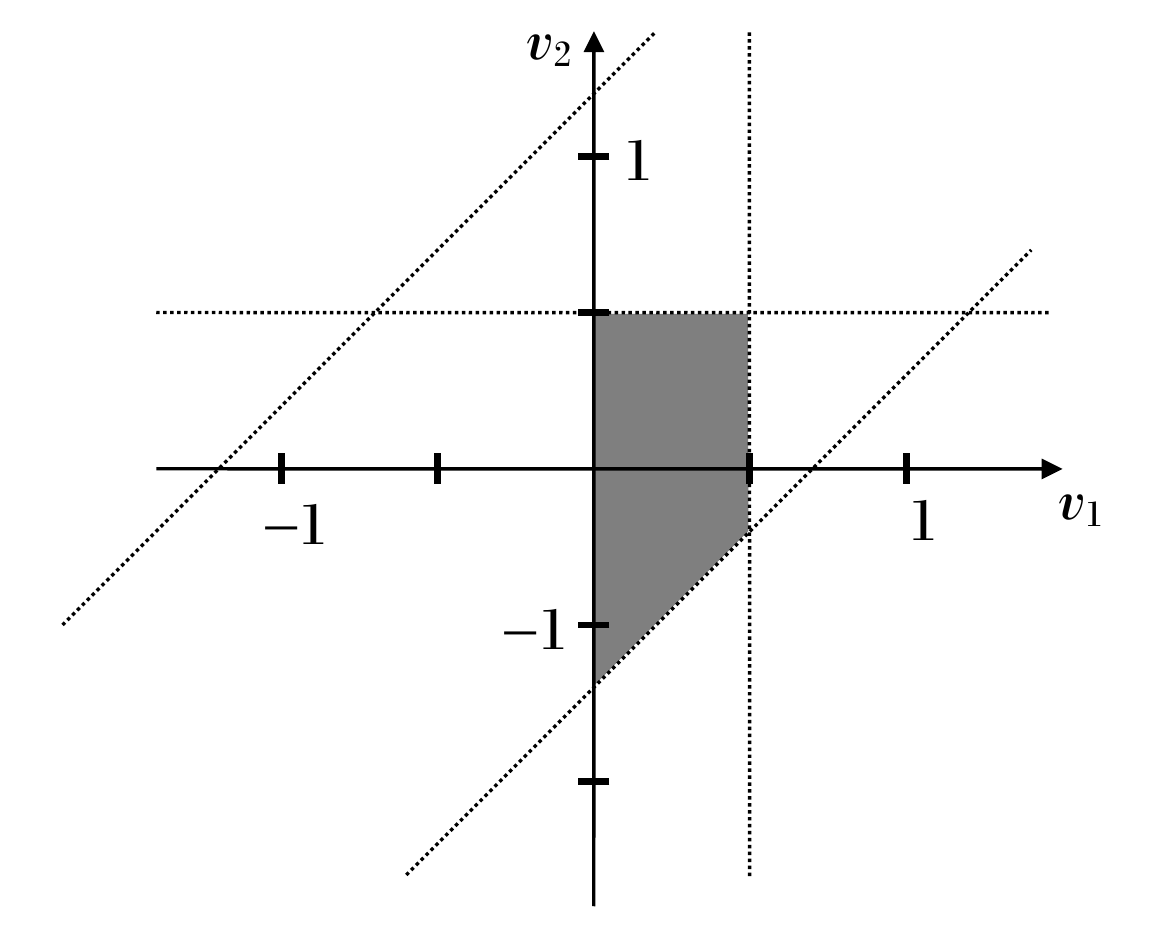}\vskip-0.2cm
\caption{\footnotesize The gray region stays for the set of $(v_1,v_2)$ given by inequalities (\ref{system2}). The dotted lines denote the sets of $(v_1,v_2)$ where equalities of (\ref{system2}) are attained (the dotted line $v_1=0$ coincides with the vertical axis and another dotted line $v_2=-1.5$ is not shown). 
} \label{figshape}
\end{figure}

\noindent Fig.~\ref{figshape} illustrates that the two constraints  from \eqref{222} corresponding to normal vectors $n_1$ and $n_4$ constitute the opposite sides of the shape 
 $\Pi(t)\cap V$. This properties persists under small perturbations of the parameters (\ref{para}). Indeed, formulas (\ref{222}) and (\ref{relation2}) imply that small perturbations of the parameters (\ref{para}) lead to small parallel displacements of the dotted lines of Fig.~\ref{figshape} (without rotations), so that the two opposite parallel sides will stay. The proof of the proposition is complete. \qed

\vskip0.2cm

\noindent In order to obtain the existence of a structurally stable family of non-stationary periodic solutions it is now remains to apply the displacement-controlled loading (\ref{relation1}) of sufficiently large amplitude. We will now use Proposition~\ref{prop3} to give an estimate for the required amplitude. In the case of a 5-spring network, formula (\ref{formulaprop3}) of Proposition~\ref{prop3}  follows from
\begin{equation}\label{fofo}
\sum_{i=1}^5\dfrac{1}{a_i}\left(c_i^+-c_i^-\right)^2<\left\|V_{basis}\bar L\right\|^2_A\cdot\left(l(t_1)-l(t_2)\right)^2.
\end{equation}
In the case of parameters (\ref{para}), formula (\ref{fofo}) reduces to 
$$
\sum_{i=1}^5 2^2<\|n_1\|^2\cdot\left(l(t_1)-l(t_2)\right)^2,
$$
where $n_1$ is given by (\ref{n1}), or simply
$$
   \dfrac{160}{3}<\left(l(t_1)-l(t_2)\right)^2.
$$
Since $\frac{160}{3}\approx 53.3,$ we introduce $l(t)$ as follows
\begin{equation}\label{loadt}
   l(t)=\left\{\begin{array}{ll} 
      t, & \ {\rm if}\ t\in[0,54],\\
      -t+54, & \ {\rm if}\ t\in[54,108],
      \end{array}\right.
\end{equation}
extended to $[0,\infty)$ by 108-periodicity.

\begin{corollary} \label{corr1} Consider the network of elastoplastic springs of Fig.~\ref{ex2newfig} with the parameters (\ref{para}).  Assume the displacement-controlled loading given by (\ref{loadt}), so that $T=108.$ Then, for any parameters $a_i,$  $c_i^-,c_i^+$, $i\in\overline{1,5}$, and any  Lipschitz-continuous functions $T$-periodic $H:[0,T]\mapsto\mathbb{R}^3,$ that are close to those in (\ref{para}), and for any Lipschitz-continuous $T$-periodic $l(t)$ close to (\ref{loadt}), the sweeping process  
(\ref{sw1})-(\ref{eq:initialElastic}) admits a structurally stable family of non-stationary $T$-periodic solutions (swept by the opposite parallel sides of Fig.~\ref{figshape}). Accordingly, the mechanical model of Fig.~\ref{ex2newfig} admits an entire family of co-existing stress distributions that evolves $T$-periodically in time. 
\end{corollary}

\section{Conclusions} In this paper we showed that sweeping processes of networks of elastoplastic springs (elastoplastic systems) inherit a designated structure that restrict possible dynamic transitions. Specifically, we gave an example of an elastoplastic system whose sweeping process admits a structurally stable family of non-stationary periodic solutions. Specifically, the structure given by the elastoplastic system locks the family of periodic solutions of the associated sweeping process, so that it persists under all such small perturbations of the sweeping process that come from small perturbations of the physical parameters of the elastoplastic system.

\sloppy

\section*{Compliance with Ethical Standards}

\noindent {\bf Conflict of Interest:} The authors have no conflict of interest.

\bibliographystyle{plain}


\end{document}